\def \version {2019--10--30}

\def \msk {\medskip}
\def \bsk {\bigskip}

\def \ex {\mathrm{ex}}
\def \rex {\mathrm{rex}}
\def \rps {K_{r+1}}

\newtheorem{theorem}{Theorem}
\def \thm {\begin{theorem}}
\def \ethm {\end{theorem}}

\newtheorem{prop}{Proposition}
\def \prp {\begin{prop}}
\def \eprp {\end{prop}}

\newtheorem{probl}{Problem}
\def \prm {\begin{probl}}
\def \eprm {\end{probl}}

\newtheorem{conjecture}{Conjecture}
\def \cnj {\begin{conjecture}}
\def \ecnj {\end{conjecture}}

\newcommand{\pft}[1]{{\bsk\noindent\bf Proof of Theorem #1}\smallskip }
\newcommand{\pfp}[1]{{\bsk\noindent\bf Proof of Proposition #1}\smallskip }
\newcommand{\pfr}[1]{{\bsk\noindent\bf Proof of Part $\bm{(#1)}$}\smallskip }

\documentclass[12pt]{article}
\usepackage{amssymb}
\usepackage{bm}

\begin{document}

\title{Regular Tur\'an numbers}
\author{Yair Caro~\thanks{Department of Mathematics,
 University of Haifa-Oranim, Tivon 36006, Israel.
  E-mail:
  {\tt yacaro@kvgeva.org.il}}
 \qquad and \qquad
  Zsolt Tuza~\thanks{Alfr\'ed R\'enyi Institute of Mathematics,
        H--1053 Budapest,
 Re\'altanoda u.~13--15, Hungary; and
  Department of Computer Science and Systems Technology,
  University of Pannonia, 8200 Veszpr\'em, Egyetem u.~10,
 Hungary.
  E-mail: {\tt tuza@dcs.uni-pannon.hu}.
 Research supported in part by the National Research,
 Development and Innovation Office -- NKFIH under the grant SNN 129364.} }
\date{\small Latest update on \version}
\maketitle

\vspace{-7ex} ~~

\begin{abstract}
 The regular Tur\'an number of a graph $F$,
  denoted by $\rex(n,F)$,
 is the largest number of edges in a regular graph $G$
 of order $n$ such that $G$ does not contain subgraphs
  isomorphic to $F$.
 Giving a partial answer to a recent problem raised by
  Gerbner et al.\ [arXiv:1909.04980]
  we prove that $\rex(n,F)$ asymptotically equals the
   (classical) Tur\'an number
    whenever the chromatic number of $F$ is at least four;
   but it is substantially different for some 3-chromatic graphs $F$
    if $n$ is odd.
\end{abstract}

\section{Introduction}

Let $F$ be a fixed `forbidden' graph. We denote
 \begin{itemize}
   \item $\ex(n,F)$ the maximum number of edges in a graph of order $n$ that
    does not contain $F$ as a subgraph --- the classical \emph{Tur\'an number};
   \item $\rex(n,F)$ the maximum number of edges in a \emph{regular} graph of order $n$ that
    does not contain $F$ as a subgraph --- the \emph{regular Tur\'an number}.
 \end{itemize}
Of course $\rex(n,F)\leq \ex(n,F)$ holds for every $F$ by definition.

The Tur\'an number of graphs is one of the most famous functions of graph theory,
 and it is well known to satisfy
  $$
    \ex(n,\rps) =
     {n\choose 2} - \sum_{i=0}^{r-1} {\lfloor \frac{n+i}{r} \rfloor \choose 2}
      \quad \mathrm{and} \quad
    \ex(n,F)=(1+o(1))\;\ex(n,K_{\chi(F)})
  $$
 for all $F$ with chromatic number $\chi(F) = r+1 \geq 3$,
 by the celebrated theorems of Tur\'an\footnote{The unique extremal graph
  for $\rps$ --- the
   \emph{Tur\'an graph}, often denoted by $T_{n,r}$ ---
  is obtained by partitioning the $n$ vertices into $r$ classes as equally as possible
   (each class has $\lfloor n/r \rfloor$ or $\lceil n/r \rceil$ vertices), and
   two vertices are adjacent if and only if they belong to distinct classes.}
  \cite{Tu41} and Erd\H os and Stone \cite{ESt46}.

The regular Tur\'an number was introduced recently by Gerbner,
 Patk\'os, Vizer, and the second author in \cite{GPTV},
 motivated by the study of singular Tur\'an numbers introduced in \cite{CT}.
It was proved in \cite{GPTV} that
 \begin{itemize}
   \item $\rex(n,K_3)$ is not a monotone function of $n$ as
     $\rex(n,K_3)=\ex(n,K_3)=n^2\!/4$ if $n$ is even,
      while $$\rex(n,K_3)\leq 2 n^2\!/5$$ if $n$ is odd,
      the latter derived from a theorem of Andr\'asfai \cite{An64},
     which we shall present in details later as well as some
      recent progress motivated by this theorem;
   \item there exists a quadratic lower bound on $\rex(n,F)$
     whenever $\chi(F)\geq 3$,
     namely $$\rex(n,F) \geq n^2 /(g + 6) - O(n)$$ where $g$ is the length of a
     shortest \emph{odd} cycle in $F$ (that is, the \emph{odd girth} of $F$);
   \item if $\chi(F)=r+1\geq 3$ and $n$ is a multiple of $r$, then
     $$\rex(n,F) = (1 + o(1))\;\ex(n,F)$$ as $n\to\infty$, by the regularity
     of the Tur\'an graph;
   \item if $F$ is a tree on $p+1$ vertices and $\ex(n,F) \leq (p - 1)n/2$,
     then $\ex(n,F) = \rex(n,F)$ for every $n$ divisible by $p$.
 \end{itemize}
Problem 4.2 of \cite{GPTV} asks to determine \ $\liminf \rex(n,F)/n^2$ \
 for non-bipartite graphs $F$.
The goal of our present note is to solve this problem for a large class of graphs $F$,
 as expressed in the following results.

\thm
 \label{t:4kr}
 Let\/ $r\geq 3$.
  \begin{itemize}
    \item[$(i)$] If\/ $n$ is a multiple of\/ $r$, then\/
      $\rex(n,\rps) = \ex(n,\rps)$.
    \item[$(ii)$] If\/ $n$ is not a multiple of\/ $r$, then\/
      $\rex(n,\rps) = \ex(n,\rps) - \Theta(n)$ as\/ $n\to\infty$.
    \item[$(iii)$] More exactly, if\/ $n=qr+s$ with\/ $1\leq s\leq r-2$,
     and at least one of\/ $r-s$ and\/ $q$ is even, then\/
     $$\rex(n,\rps) = \ex(n,\rps) - \frac{(r-s) q}{2}  .$$
    \item[$(iv)$] If\/ $F$ is any graph with\/ $\chi(F)\geq 4$,
     then\/ $\rex(n,F) = (1 - o(1))\;\ex(n,F)$.
    \item[$(v)$] For\/ $n=3k+s$ with\/ $s=0,1,2$
     we have\/ $\rex(n,K_4) = kn = n\cdot \lfloor \frac{n}{3} \rfloor$.
  \end{itemize}
\ethm

\thm Let\/ $F$ be a\/ $3$-chromatic graph.
 \label{t:c3}
  \begin{itemize}
    \item[$(i)$] If\/ $n$ is even, then\/
      $\rex(n,F) = (1 - o(1))\;\ex(n,F) = n^2\!/4 + o(n^2)$;
     moreover,\/ $\rex(n,K_3) = \rex(n,K_4-e) = n^2\!/4$.
    \item[$(ii)$] If\/ $n$ is odd, and\/ $F=K_3$ or\/ $F=K_4-e$
      or\/ $F$ is a unicyclic graph with\/ $C_3$ as its cycle,
      then\/ $\rex(n,F) = n^2\!/5 - O(n)$.
    \item[$(iii)$] If\/ $F=K_3$ and\/ $n=5k+s$ is odd,
      then\/ $\rex(n,F) = kn = n\cdot \lfloor \frac{n}{5} \rfloor$.
  \end{itemize}
\ethm

\thm
 \label{t:cg}If\/ $\chi(F)=3$, and\/ $F$ has odd girth\/ $g$, then\/
     $$\rex(n,F) \geq n^2\!/(g+2) - O(n) $$
   for every odd\/ $n$.
  Moreover, if\/ $F=C_5$
      or\/ $F$ is a unicyclic graph with\/ $C_5$ as its cycle,
      then the asymptotic equality\/ $$\rex(n,F) = n^2\!/7 - O(n)$$
       is valid.
  More precisely, if\/ $n=7k+s$ is odd, and\/ $n$ is sufficiently large
    with respect to\/ $F$,
      then\/ $\rex(n,F) = kn = n\cdot \lfloor \frac{n}{7} \rfloor$.
\ethm
 
It should be noted that almost nothing is known about the behavior of
 $\rex(n,F)$ in case of bipartite graphs $F$, apart from a couple of
 remarks in the concluding section of \cite{GPTV}.
Moreover, the following problem remains widely open for 3-chromatic graphs.

\prm
 \label{c:g}
 Determine\/ $\rex(n,F)$, or its asymptotic growth as\/ $n\to\infty$,
  for graphs\/ $F$ with\/ $\chi(F)=3$ for odd\/ $n$.
\eprm

As a particular case, we expect that the exact results on $K_3$ and $C_5$
 extend to a unified formula for every odd cycle.

\cnj
 \label{c:cyc}
 If\/ $g\geq 3$ and\/ $n=k\, (g+2) + s$ with\/
  $0\leq s\leq g+1$, and both\/ $g$\/ and $n$\/ are odd, then\/
  $\rex(n,C_g) = kn = n\cdot \lfloor \frac{n}{g+2} \rfloor$.
\ecnj

The unicyclic extensions given in Theorems \ref{t:c3} and \ref{t:cg}
 are consequences of the following principle.
It has an analogous implication also for those graphs whose unique
 non-trivial block is $K_4-e$.

\prp
 \label{p:pend}
 If the growth of\/ $\rex(n,F)$ is superlinear in\/ $n$, and\/ $F^+$
  is a graph obtained from\/ $F$ by inserting a pendant vertex, then\/
   $\rex(n,F^+) = \rex(n,F)$ for every sufficiently large\/ $n$.
\eprp

At the end of this introduction let us recall the full statement of
 And\-r\'as\-fai's theorem, which plays an essential role in the current context.

\thm [\cite{An64}]
 \label{t:And}
 If\/ $G$ is a triangle-free graph on\/ $n$ vertices and with
  minimum degree\/ $\delta(G) >2n/5$, then\/ $G$ is bipartite.
\ethm

This theorem is the source of motivation for recent research which is
 also related and relevant to our Theorem \ref{t:cg} and
  Conjecture \ref{c:g}; see \cite{ABGKM,EbSch,LetSny,MesSch}.
We explicitly quote the following very useful generalization,
 proved by Andr\'asfai, Erd\H os, and T. S\'os, as read out from the
  combination of their Theorem 1.1 and Remark 1.6.

\thm [\cite{AES}]
 \label{t:AES}

For each odd integer\/ $k\geq 5$ and each integer\/ $n\geq k$,
 if\/ $G$ is a simple\/ $n$-vertex graph with no odd cycles of length
  less than\/ $k$ and with minimum degree $\delta(G)>2n/k$,
   then\/ $G$ is bipartite.
\ethm

\section{Forbidden graphs with chromatic number at least 4}

In this section we prove Theorem \ref{t:4kr}.
We assume throughout that $n=qr+s$, where $q$ is an integer and
 $0\leq s\leq r-1$ holds.

\pfr{i}

Whenever $n$ is a multiple of $r$, the equality
 $\rex(n,\rps) = \ex(n,\rps)$ is clear because the $r$-partite
 Tur\'an graph $T_{n,r}$ is regular in all such cases.

\pfr{ii}

First we argue that $\ex(n,\rps) - \rex(n,\rps)$
  is at least a linear function of $n$
 if $r$ does not divide $n$.
We see from Tur\'an's theorem that the degree of regularity cannot exceed
 $(1-1/r)\,n$ if the graph is $\rps$-free.
For $n=qr+s$ with $0<s<r$ it means that the degrees are at most $n-(n+r-s)/r$.
This value is the degree of vertices in the $s$ larger classes of the Tur\'an graph;
 but the $r-s$ smaller classes consist of vertices of degree $n-(n+r-s)/r+1$.
Hence the degree sum in a regular graph is smaller by at least $\frac{r-s}{r}\, n - O(1)$.

\msk

Next, we construct $r$-chromatic (hence also $\rps$-free)
 regular graphs to show that the difference
 between $\ex(n,\rps)$ and $\rex(n,\rps)$ is at most $O(n)$.
For $n=qr+s$ with $0<s<r$ the Tur\'an graph has $s$ classes of size $q+1$
 and $r-s$ classes of size $q$.
Putting this in another way, the vertices in $s$ classes have degree $n-q-1$, and
 in $r-s$ classes have degree $n-q$.

We are going to delete $O(n)$ edges and obtain a regular graph.
For this purpose we shall use Dirac's theorem \cite{Di52}
 as a lemma, which states that if
 the degree of every vertex in a graph $H$ is at least half of the number of vertices,
 then $H$ has a Hamiltonian cycle.

Assume first that both $s$ and $r-s$ are at least 2.
Since the degree sum is even, the total size of the odd-degree classes together
 is also even, and the subgraph induced by them is Hamiltonian.
Hence the union of these classes admits a 1-factor, which we remove.
If this is the larger degree, then we are done.
Otherwise the subgraph induced by the larger degrees also has a Hamiltonian cycle,
 whose removal leaves a regular graph of degree $n-q-2$.

Assume next that $s=1$; i.e., only one class contains vertices of low degree.
If the total size of the $r-1$ high-degree classes is even, we remove a 1-factor
 from the subgraph induced by them, and we are done.
On the other hand, if their total size is odd, their degree must be even.
That is, both $q$ and $r-1$ are odd; in particular, $r-1\geq 3$ holds.
In this situation our plan is to delete two edges from each such vertex,
 and one edge from each vertex of the low-degree class.

We begin with the single class, which has odd size.
We omit one edge from each of its vertices --- mutually disjoint edges ---
 in such a way that the other ends of those $q+1$ edges are distributed
 as equally as possible among the $r-1\geq 3$ high-degree classes.
Then Dirac's condition is valid for the high-degree ends of the omitted
 matching and also for the subgraph from which no edges have been
 omitted so far.
Hence both parts are Hamiltonian.
Since $q+1$ is even, we can omit a perfect matching from the former,
 and a Hamiltonian cycle from the latter, thus obtaining a regular graph.

Finally, consider the case of $s=r-1$; i.e., only one class contains vertices
 of high degree.
We must remove edges from the high-degree class, which means that also some low-degree
 vertices will decrease their degree.
Similarly to the previous case, here again, the parity of classes will matter.

If the high-degree class has even size $q$, we decrease its degrees by 3,
 distributing the neighbors equally among the low-degree classes, and inside the
 neighbors we delete a 1-factor.
In the rest of low-degree vertices we delete a Hamiltonian cycle.

Suppose that the high-degree class has odd size $q$.
Then the number of low-degree vertices is even because each such class has
 even size $q+1$.
In particular, $n$ is odd.
We now delete two edges from each high-degree vertex, this decreases $2q$ of the
 low degrees by 1.
From the other $n-3q$ (i.e., even number of) vertices we delete a 1-factor
 which exists since also this subgraph is Hamiltonian.
This modification yields a regular graph, and completes the proof.

\pfr{iii}

Recall that the number of classes of high-degree vertices in the Tur\'an graph
 is $r-s$, and these classes have cardinality $q$ each.
Note further that their union induces a Hamiltonian subgraph whenever
 $s\leq r-2$.
Under the assumption that at least one of $r-s$ and $q$ is even, the length
 $(r-s)\, q$ of a corresponding Hamiltonian cycle is even, hence
  contains a perfect matching, say $M$.
Removing $M$ from $T_{n,r}$ we obtain a $\rps$-free regular graph,
 and the degree is largest possible, according to the first part of the
  proof of $(ii)$ as given above.

\newpage

\pfr{iv}

The $r$-colorable construction given above for $(ii)$ proves that
 $\rex(n,F)$ is at least $\ex(n,\rps) - O(n)$ if $\chi(F)=r+1$.
On the other hand, as mentioned already,
 $\ex(n,F) = (1 + o(1))\;\ex(n,\rps)$ is valid by the
 Erd\H os--Stone theorem, implying that $\rex(n,F)$ cannot be larger.

\vbox{
\pfr{v}

The following list is a summary of optimal constructions according to $n$~(mod 3).
}

\begin{itemize}
  \item $n=3k$ --- the complete 3-partite graph with equal classes,
    i.e.\ $T_{n,3}$, is regular.
  \item $n=3k+1$ --- here $T_{n,3}$ has vertex classes of respective sizes $k,k,k+1$;
    it can be made regular by removing a perfect matching between the two classes of size $k$.
  \item $n=3k+2$ --- here $T_{n,3}$ has vertex classes of respective sizes $k,k+1,k+1$;
    it can be made regular by removing a matching of $k$ edges between the class
     of size $k$ and each class of size $k+1$ (hence removing $kP_3$ from $T_{n,3}$),
      moreover deleting the edge that joins the two
     vertices whose degree has not been decreased by the removal of the two matchings.
    This is optimal because removing an edge from a vertex of high degree decreases
     the degree of a vertex of low degree, too.
\end{itemize}

\section{3-chromatic forbidden graphs}

Let us begin this section with the proof of Proposition \ref{p:pend},
 as it is applicable for Theorems \ref{t:c3} and \ref{t:cg} as well.

\pfp{\ref{p:pend}}

 Certainly we have $\rex(n,F^+) \geq \rex(n,F)$.
Suppose that $G$ is a regular graph of order $n$,
 which is extremal for $F^+$.
If $G$ is $F$-free, then the reverse inequality $\rex(n,F^+) \leq \rex(n,F)$
 also holds and the assertion follows immediately.
On the other hand, if $F\subset G$ but $F^+\not\subset G$,
 the degree of regularity in $G$ must be smaller than $|V(F^+)|$, for
 otherwise it would be possible to extend $F$ to $F^+$ in $G$.
This fact puts a $O(n)$ upper bound on $|E(G)|$, contradicting the
 superlinear growth of $\rex(n,F)$.
Thus the largest $F^+$-free regular graphs are $F$-free, too, as $n$ gets large.

\pft{\ref{t:cg}, general lower bound}

We construct a graph of order $n$ and odd girth $g+2$, hence it will not contain $F$
 as a subgraph.
Let us write $n$ in the form $n=(g+2)\cdot a + 2b$ where $b$ is an integer in the
 range $0\leq b\leq g+1$.
We start with a blow-up of $C_{g+2}=v_1 v_2\dots v_{g+2}$ by substituting
 independent sets $A_1,A_2,\dots,A_{g+2}$ into its vertices, completely joining
 $A_i$ with $A_{i+1}$ for $i=1,\dots,g+2$ (where $A_{g+3}:=A_1$).
We let $|A_1|=|A_2|=a+b$, and $|A_i|=a$ for all $3\leq i\leq g+2$.
The degree of a vertex $v$ in this graph is $2a+b$ if
 $v\in A_1\cup A_2\cup A_3\cup A_{g+2}$, and it is $2a$ otherwise.
The graph is regular if $b=0$, and it will be made regular by the removal of
 $$
   2ab+b^2 = b \left( \frac{2n-4b}{g+2} + b \right) =
  \frac{b\cdot ( 2n + (g-2)b )}{g+2} \leq 2n - \frac{2n-(g-2)(g+1)^2}{g+2}
 $$
  edges otherwise.

Between $A_{g+2}$ and $A_1$ we remove a bipartite graph $H$ such that all vertices
 of $H$ in $A_{g+2}$ have degree $b$, and all degrees in $A_1$ are
 $\lfloor \frac{ab}{a+b} \rfloor$ or $\lceil \frac{ab}{a+b} \rceil$.
Such $H$ clearly exists.
We also remove a bipartite graph isomorphic to $H$ between $A_2$ and $A_3$, such that
 the degree-$a$ vertices are in $A_3$.

If $ab$ is a multiple of $a+b$, then the current vertex degrees in $A_1\cup A_2$ are
 $2a+b - \frac{ab}{a+b} = 2a + \frac{b^2}{a+b}$.
Then removing a regular bipartite graph of degree $\frac{b^2}{a+b}$, hence with
 $b^2$ edges, yields a $(2a)$-regular graph of order $n$.

Otherwise, if $ab$ is not divisible by $a+b$, we first remove a perfect matching
 between the vertices of degree $2a+b - \lfloor \frac{ab}{a+b} \rfloor$ in $A_1\cup A_2$.
After that, the bipartite graph induced by $A_1\cup A_2$ is regular, hence
 deleting a regular subgraph of degree $b - \lceil \frac{ab}{a+b} \rceil$ from it,
 we obtain a $(2a)$-regular graph of order $n$.

\pft{\ref{t:cg} for asymptotics of $\bm{C_5}$ and unicyclic graphs}

 We have to prove that $\rex(n,F) \leq n^2\!/7$, if $n$ is not very small.
By Proposition \ref{p:pend} it is enough to deal with the case of
 $F=C_5$.
Let $G$ be a $C_5$-free regular graph with $\rex(n,C_5)$ edges.
 If $G$ is triangle-free, the proof is done by Theorem \ref{t:AES},
 because the odd girth cannot be exactly 5.
On the other hand it was proved in \cite[Lemma 33]{LetSny} that if $G$ contains
 a triangle and the degrees are greater than $(1/6+\epsilon)\,n$,
 for any $\epsilon>0$ and sufficiently large $n$, then $G$ also has
 a $C_5$.
Hence in $C_5$-free graphs with triangles we cannot have more than
 $n^2\!/12 + o(n^2)$ edges.

We postpone the proof of the exact formula for $\rex(n,C_5)$ to the
 end of this paper,
 due to its similarity to the argument concerning $\rex(n,C_3)$.

\pft{\ref{t:c3}, Parts $\bm{(i)}$ and $\bm{(ii)}$}

If $F$ is 3-chromatic and $n$ is even, the lower bound of $n^2\!/4$
 is shown by the complete bipartite graph $K_{n/2,n/2}$, while an
 asymptotic upper bound follows by the Erd\H os--Stone theorem as
  $$
    \rex(n,F) \leq \ex(n,F) = (1 - o(1))\;\ex(n,K_3) = n^2\!/4 + o(n^2) .
  $$
The tight results $\rex(n,K_3) = \rex(n,K_4-e) = n^2\!/4$ follow from the
 facts that the Tur\'an number of $K_3$ and also of $K_4-e$ is $n^2\!/4$.

 Assume that $n$ is odd.
The lower bound $\rex(n,K_3) \geq n^2\!/5 - O(n)$ is a particular case
 of the previous construction, putting $g=3$.
Moreover, in triangle-free regular graphs we cannot have more than
 $n^2\!/5$ edges, due to Theorem \ref{t:And} and by the fact that
 every regular bipartite graph has an even order.
This already settles the case of $K_3$ (and also of the unicyclic graphs
 having a 3-cycle).
For $K_4-e$ assume that $x,y,z$ induce $K_3$.
If this triangle cannot be extended to $K_4-e$, then the degree of
 $x,y,z$ is at most $2+(n-3)/3 = n/3 +1$, thus
 by the condition of regularity the number of edges is
 at most $n^2\!/6 + n/2$, which is much less than $n^2\!/5$ if $n$ is large.

\pft{\ref{t:c3}, Part $\bm{(iii)}$}

Let $G$ be a triangle-free regular graph on $n$ vertices, with
 $|E(G)|=\rex(n,K_3)$.
Assume that $n = 5k +s$, where $s = 0,1,2,3,4$ and $n$ is odd.
As a consequence, $k+s$ is odd as well; however, the degree $d$
 of regularity must be even.
From Theorem \ref{t:And} we also know that
 $d\leq \lfloor 2n/5 \rfloor = 2k + \lfloor 2s/5 \rfloor$.
A more careful look verifies that the last term $\lfloor 2s/5 \rfloor$
 will vanish with respect to $d$, and in every possible case we have
  $$
    d \leq 2k = 2(n-s)/5 .
  $$
 This is clear for $s=0,1,2$ due to the floor function.
But also for $s=3,4$ we have that $\lfloor 2s/5 \rfloor = 1$ is an
 odd number, whereas $d$ must be even, hence $d$ cannot exceed the
 largest even integer under $2k+1$, which is actually $2k$.
Thus, $|E(G)|\leq kn$.

It remains to show that for every odd $n = 5k +s$ there exists a
 $K_3$-free graph of order $n$ which is $2k$-regular.
The general principle of the construction is to substitute
 independent sets $A_1,\dots,A_5$ into the vertices of $C_5$,
 where each edge of $C_5$ becomes a complete bipartite graph
  between the corresponding two sets $A_i,A_{i+1}$ cyclically;
 and then delete some edges so that a regular graph is obtained.
We are going to describe these constructions for each $s$
 one by one, specifying the sequences $(|A_1|,\dots,|A_5|)$ as follows.
  \begin{itemize}
    \item $s=0$ :\quad $(|A_1|,|A_2|,|A_3|,|A_4|,|A_5|) = (k,k,k,k,k)$

      This graph is $2k$-regular.

    \item $s=1$ :\quad $(|A_1|,|A_2|,|A_3|,|A_4|,|A_5|) = (k+1,k+1,k,k-1,k)$

      This graph becomes $2k$-regular after the deletion of a perfect
       matching from the induced subgraph $G[A_1\cup A_2]$.

    \item $s=2$ :\quad $(|A_1|,|A_2|,|A_3|,|A_4|,|A_5|) = (k+1,k+1,k,k,k)$

      This graph becomes $2k$-regular after the deletion of a
       matching of size $k$ from $G[A_1\cup A_5]$ and from $G[A_2\cup A_3]$,
       moreover the edge between the two unmatched vertices of $A_1\cup A_2$.

    \item $s=3$ :\quad $(|A_1|,|A_2|,|A_3|,|A_4|,|A_5|) = (k+1,k+1,k+1,k,k)$

      This graph becomes $2k$-regular after the deletion of a perfect
       matching from each of the induced subgraphs
        $G[A_1\cup A_2]$, $G[A_2\cup A_3]$, $G[A_4\cup A_5]$.

    \item $s=4$ :\quad $(|A_1|,|A_2|,|A_3|,|A_4|,|A_5|) = (k+2,k+2,k,k,k)$

      This graph can be made $2k$-regular in the following way.
      Specify two vertices $a_1',a_1''\in A_1$ and $a_2',a_2''\in A_2$;
       set $A_1'=A_1\setminus \{a_1',a_1''\}$ and $A_2'=A_2\setminus \{a_2',a_2''\}$.
      Delete a 2-factor from $G[A_1'\cup A_5]$ and from $G[A_2'\cup A_3]$;
       and delete the edges of the 4-cycle $a_1'a_2'a_1''a_2''$.
  \end{itemize}

\pft{\ref{t:cg} for exact $\bm{\rex(n,C_5)}$}

Since the proof is very similar to that of the exact formula for $\rex(n,K_3)$,
 we give a more concise description here.
  Let now $n=7k+s$, where $0\leq s\leq 6$.
We have already seen that the degree $d$ of regularity satisfies
 $d \leq \lfloor 2n/7 \rfloor = 2k + \lfloor 2s/7 \rfloor \leq 2k+1$;
 and $d$ must be even, thus $d\leq 2k$.
It remains to give suitable substitutions of sets $A_1,\dots,A_7$ into the
 vertices of $C_7$ in such a way that the graphs can be made $2k$-regular
 by the deletion of some edges.
Below we define a sequence $|A_1|,|A_2|,|A_3|,|A_4|,|A_5|,|A_6|,|A_7|$
 similar to the case of $C_5$, now for each $s=0,1,\dots,6$.
\begin{itemize}
 \item $s=0$ :\quad $k,k,k,k,k,k,k$

   Nothing to delete.
 \item $s=1$ :\quad $k+1,k,k,k+1,k,k-1,k$

   Delete a 1-factor from $G[A_2\cup A_3]$.
 \item $s=2$ :\quad $k+1,k+1,k,k,k,k,k$

   Delete a matching of size $k$ from $G[A_1\cup A_7]$ and from $G[A_2\cup A_3]$,
    and the edge between the two unmatched vertices of $A_1\cup A_2$.
 \item $s=3$ :\quad $k+1,k+1,k,k,k+1,k,k$

   Delete a 1-factor from $G[A_1\cup A_2]$, from $G[A_3\cup A_4]$, and
    from $G[A_6\cup A_7]$.
 \item $s=4$ :\quad $k+2,k+2,k,k,k,k,k$

   Delete the edges of a $H\cong C_4$ in $A_1\cup A_2$, and a 2-factor from
    $G[(A_1\cup A_7)\setminus V(H)]$ and from $G[(A_2\cup A_3)\setminus V(H)]$.
 \item $s=5$ :\quad $k+2,k+2,k,k,k+1,k,k$

   Delete a $C_4$ from $G[A_1\cup A_2]$, and a matching of size $k$ from each
    consecutive pair of $A_i,A_{i+1}$ along the cycle (also including $A_7,A_1$
    as cyclically consecutive), except $G[A_4\cup A_5]$ and $G[A_5\cup A_6]$.
 \item $s=6$ :\quad $k+2,k+2,k,k,k+2,k,k$

   Delete a 2-factor from $G[A_1\cup A_2]$, from $G[A_3\cup A_4]$, and
    from $G[A_6\cup A_7]$.
\end{itemize}
After the deletions, all graphs are $2k$-regular, completing the proof of the
 theorem.
\bsk

\paragraph{Acknowledgements.}

We thank D\'aniel Gerbner, Bal\'azs Patk\'os, and M\'at\'e Vizer for discussions 
 during the preparation of this manuscript.

\end{document}